# A Bare Bones Grey Wolf Optimizer for Global Numerical Optimization

Haoxin Wang
*National Key Laboratory of Power Systems in Shenzhen*
*Shenzhen International Graduate School, Tsinghua University*
Shenzhen, China
whx980425ee@163.com

Libao Shi
*National Key Laboratory of Power Systems in Shenzhen*
*Shenzhen International Graduate School, Tsinghua University*
Shenzhen, China
shilb@sz.tsinghua.edu.cn

*Abstract*—In order to better understand and analyze the currently widely used population-based metaheuristic optimization algorithms, , this paper proposes a novel computational intelligence algorithm called bare bones grey wolf optimizer (BBGWO) inspired by a bare bones mechanism. In the BBGWO, the complex updating mechanism of solutions is replaced by a random vector that obeys a normal distribution, whose mean and variance are derived by theoretically studying the probability distribution of the new solution of the original GWO. The corresponding theoretical analysis and simulation results verify the good optimization performance of the proposed BBGWO algorithm .

## I. Introduction

The performance of traditional mathematical programming methods (mostly represented by conjugate gradient method and Newton method) highly depends on the structure of the optimization problem to be solved, and might encounter solution bottlenecks when solving multimodal or highly non-convex optimization problems. In order to deal with the increasingly complex engineering optimization problems effectively, in the past fifty years, a large number of population-based metaheuristic optimization algorithms have gained rapid development and widespread attention, such as genetic algorithm (GA) [1], particle swarm optimization (PSO) [2], ant colony optimization (ACO) [3], etc. Most of these algorithms are based on a metaphor of some natural or man-made process [4], and hold solution searching mechanisms that are hardly dependent on the structure and characteristics of optimization problems to be solved. These features guarantee their extremely attractive advantages in solving highly complex or multimodal optimization problems. However, since the solution updating mechanisms of these algorithms are usually designed by mimicking some natural or man-made process, the expressions of new solutions might be non-intuitive and hard to understand, making it difficult to conduct in-depth study on the mechanism of these algorithms (for instance, the trajectory of a single individual, the topological structure of the entire population, etc.).

To better understand and analyze these population-based metaheuristic algorithms, James Kennedy [5] proposed a bare bones mechanism for simplifying these algorithms, and applied this mechanism to the original PSO to obtain a modified PSO called bare bones PSO (BBPSO). In the BBPSO, the velocity vector of an individual was eliminated, and the solution updating equation was revised as a random vector that obeys a normal distribution (see chapter 2.1 for more details). Theoretical analysis [6] and simulation results showed that the BBPSO could highly imitate the optimization mechanism of the original PSO. In addition, owing to the intuitive solution updating strategy and the excellent performance (not worse than the original PSO algorithm in most cases [7]), the BBPSO has been widely used to solve engineering optimization problems independently as well, including unsupervised feature selection [8], data clustering [9], environmental/ economic dispatch [10], etc. Since the proposing of the BBPSO, the bare bones mechanism has been applied to many other metaheuristic algorithms for their simplification and modification as well, such as differential evolution (DE) [11], fireworks algorithm (FA) [12], brain-storm optimization (BSO) [13], etc.

As one of the most typical population-based metaheuristic optimization algorithms, the grey wolf optimizer (GWO) [14] proposed by Seyedali Mirjalili has become a powerful tool for solving engineering optimization problems like economic dispatch [15], image segmentation [16], clustering [17], etc. due to its outstanding performance [18] and fewer parameters. However, compared with its excellent ability based on simulation analysis, few studies have been conducted on its theoretical analysis.

In this paper, we apply the bare bones mechanism to the original GWO to obtain a novel algorithm called bare bones GWO (BBGWO), in which the complex solution updating mechanism is replaced by a random vector that obeys a normal distribution, whose mean and variance are derived by calculating those of the new solution of the original GWO. Theoretical analysis and simulation results show that the mechanism of the proposed BBGWO algorithm is highly similar to that of the original GWO, which means the proposed BBGWO holds both theoretical value of serving as a fundamental tool in studying the mechanism of the GWO and practical significance of being applied to solve engineering optimization problems independently.

## II. Preliminaries

For a better illustration of the proposed BBGWO algorithm, in this section, we first briefly introduce the bare bones mechanism and the BBPSO proposed in [5], and then briefly introduce the original GWO algorithm.

### A. A brief introduction of bare bones mechanism and BBPSO

Here, we consider the *D*-dimension minimization problem with objective function $f(\mathbf{x})$. In the PSO, for a population with *N* individuals, each individual is defined with two vectors: 1) position vector $\mathbf{x_i} = [x_{i1}, ..., x_{iD}] = [x_{ij}]$ (*i* is an index of a certain individual of the population, *j* is an index of a certain dimension of the independent variable vector, vectors are shown in bold and components of vectors are shown in italics, similarly hereinafter); 2) and velocity vector $\mathbf{v_i} = [v_{ij}]$. The position vectors represent solutions of $f(\mathbf{x})$, and the velocity

vectors are intermediate variables for deriving position vectors, which are calculated as:

$$v_{ij}(t+1) = v_{ij}(t) + C_1 r_1 \left(p_{ij}^{(best)}(t) - x_{ij}(t)\right)$$
$$+ C_2 r_2 \left(p_{gj}^{(best)}(t) - x_{ij}(t)\right) \quad (1)$$

$$x_{ij}(t+1) = x_{ij}(t) + v_{ij}(t+1) \quad (2)$$

where $t$ denotes the iterations of the algorithm, $\mathbf{p_i}^{(best)} = \left[p_{ij}^{(best)}\right]$ is the best known position of particle $i$:

$$\mathbf{p_i}^{(best)}(t) = \mathbf{x_i}(t_0), t_0 = \arg\min_t f(\mathbf{x_i}(t)) \quad (3)$$

and $\mathbf{p_g}^{(best)} = \left[p_{gj}^{(best)}\right]$ is the best known position of the entire swarm:

$$\mathbf{p_g}^{(best)}(t) = \mathbf{p_{i_0}}^{(best)}(t), i_0 = \arg\min_i f\left(\mathbf{p_i}^{(best)}(t)\right) \quad (4)$$

$r_1, r_2 \sim U[0,1]$, and $C_1, C_2$ are coefficients.

Obviously, the position and velocity vectors are mixed together in the solution updating equation, causing the inconvenience of PSO theoretical analysis. In view of this, the bare bones mechanism was applied to the original PSO to revise the solution updating equation. The corresponding expressions are given as follows:

$$x_{ij}(t+1) \sim N(\mu, \sigma^2) \quad (5)$$

$$\mu = \frac{p_{ij}^{(best)}(t) + p_{gj}^{(best)}(t)}{2}, \sigma = \left|p_{gj}^{(best)}(t) - p_{ij}^{(best)}(t)\right| \quad (6)$$

where $N(\mu, \sigma^2)$ denotes the normal distribution with mean $\mu$ and standard deviation $\sigma$. This variant of PSO is called bare bones PSO (BBPSO). Obviously, the bare bones mechanism eliminates the velocity vectors, and converts the original solution updating equation into random vectors that obey normal distribution.

The simulation results show that for most cases the performance of BBPSO is similar to that of the original PSO, which means the BBPSO itself can be applied to solve optimization problems independently. Moreover, the theoretical analysis and simulation results show that the mechanism of BBPSO is highly similar to that of PSO, implying that the theoretical analysis of PSO can be carried out by studying the mechanism of BBPSO. From this perspective, the advantages of BBPSO can be summarized in following twofold: 1) the elimination of velocity vector highly simplifies the original second-order dynamical system, which has been widely accepted as a model of a single solution of the original PSO [19]; 2) the probability distribution of the new solution is directly given, which may serve as a fundamental tool for analyzing the mechanism of the algorithm.

*B. A Brief Introduction of GWO*

As one of the most typical population-based metaheuristic optimization algorithms, the GWO finds the optimal solution by mimicking the predation process of a group of grey wolves. The individuals update their positions by the following equations [14]:

$$\mathbf{A} = 2\mathbf{a} \cdot \mathbf{r_1} - \mathbf{a} \quad (7)$$

$$\mathbf{C} = 2\mathbf{r_2} \quad (8)$$

$$\begin{cases} \mathbf{D_1} = |\mathbf{C_1} \cdot \mathbf{p_1}(t) - \mathbf{x_i}(t)| \\ \mathbf{D_2} = |\mathbf{C_2} \cdot \mathbf{p_2}(t) - \mathbf{x_i}(t)| \\ \mathbf{D_3} = |\mathbf{C_3} \cdot \mathbf{p_3}(t) - \mathbf{x_i}(t)| \end{cases} \quad (9)$$

$$\begin{cases} \mathbf{x'_1}(t) = \mathbf{p_1}(t) - \mathbf{A_1} \cdot \mathbf{D_\alpha} \\ \mathbf{x'_2}(t) = \mathbf{p_2}(t) - \mathbf{A_2} \cdot \mathbf{D_\beta} \\ \mathbf{x'_3}(t) = \mathbf{p_3}(t) - \mathbf{A_3} \cdot \mathbf{D_\delta} \end{cases} \quad (10)$$

$$\mathbf{x_i}(t+1) = \frac{1}{3}\sum_{k=1}^{3} \mathbf{x_k}(t) \quad (11)$$

where $\mathbf{p_k} = [p_{kj}]$ ($k = 1,2,3$ denotes the index of the best three individuals, similarly hereinafter) denote the position vector of the best individual, second best individual and third best individual respectively, $\mathbf{D_k}, \mathbf{x'_k} = [x'_{kj}]$ are intermediate variables, $\mathbf{A} = [a, ..., a]$ is the step parameter vector, and $a$ is linearly decreased from 2 to 0 during of iteration:

$$a = 2\left(1 - \frac{t}{T}\right) \quad (12)$$

where $T$ denotes the maximum iterations. $\mathbf{A} = [A_j], \mathbf{C} = [C_j]$ are parameter vectors, and $\mathbf{r_1}, \mathbf{r_2}$ are random vectors that obey $U[0,1]$.

### III. BBGWO

In this section, we apply the bare bones mechanism to the original GWO algorithm to propose the BBGWO algorithm. More specifically, we theoretically study some major characteristics of the probability distribution of the new solution, and the corresponding simulations are carried out to verify the validity of the proposed method.

*A. A Theoretical Analysis of GWO*

Equations (7) - (11) can be rewritten as:

$$x_{ij}(t+1) = \frac{1}{3}\sum_{k=1}^{3} x'_{kj}(t)$$
$$= \frac{1}{3}\sum_{k=1}^{3} \left[p_{kj}(t) + A_{kj}\left|C_{kj}p_{kj}(t) - x_{ij}(t)\right|\right] \quad (13)$$

where $A_{kj}$ i.i.d. $\sim U[-a, a]$, $C_{kj}$ i.i.d. $\sim U[0,2]$.

For convenience, we simplify $x'_{kj}(t), p_{kj}(t), x_{ij}(t)$ as $x'_{kj}, p_{kj}, x_{ij}$, respectively. Regarding $\mathbf{p_1}, \mathbf{p_2}, \mathbf{p_3}, \mathbf{x_i}$ as constant vectors, then $x'_{kj}$ and $x_{ij}(t+1)$ can be considered as the functions of random variables $A_{kj}, C_{kj}$. Let the probability density function (PDF) of $x'_{kj}, x_{ij}(t+1)$ be $g_{kj}(u), h_{ij}(u)$, respectively, then we first derive the expression of $g_{kj}(u)$. The details are given as follows:

**Theorem 1**: Let

$$m_{kj} = a(-|p_{kj}| + |x_{ij} - p_{kj}|)$$
$$n_{kj} = a(|p_{kj}| + |x_{ij} - p_{kj}|) \quad (14)$$

then, when $m_{kj} \leq 0$,

$$g_{kj}(u) = \begin{cases} \frac{1}{n_{kj} - m_{kj}} \ln \frac{\sqrt{-m_{kj}n_{kj}}}{|u - p_{kj}|}, |u - p_{kj}| < -m_{kj} \\ \frac{1}{2(n_{\alpha i} - m_{\alpha i})} \ln \frac{n_{kj}}{|u - p_{kj}|}, -m_{kj} \leq |u - p_{kj}| < n_{kj} \\ 0, |u - p_{kj}| \geq n_{kj} \end{cases} \quad (15)$$

when $m_{\alpha i} > 0$,

$$g_{kj}(u) = \begin{cases} \frac{1}{2(n_{kj}-m_{kj})}\ln\frac{n_{kj}}{m_{kj}}, |u-p_{kj}|<m_{kj} \\ \frac{1}{2(n_{kj}-m_{kj})}\ln\frac{n_{kj}}{|u-p_{kj}|}, m_{kj}\leq|u-p_{kj}|<n_{kj} \\ 0, |u-p_{kj}|\geq n_{kj} \end{cases} \quad (16)$$

**Proof**: Mark $g_{kj}(u)$ as $g_{\alpha i}(u, x_{ij}, p_{kj})$, then according to (13), $g_{kj}(u, x_{ij}, p_{kj}) = g_{kj}(u, -x_{ij}, -p_{kj}) = g_{kj}(u, 2p_{kj} - x_{ij}, p_{kj}) = g_{kj}(2p_{kj} - u, x_{ij}, p_{kj})$, thus we only need to consider the condition $p_{kj} \geq 0, x_{ij} \geq p_{kj}$ and $u \geq p_{kj}$.

Let $D = \{(A_{kj}, C_{kj})|p_{kj} + A_{kj}|C_{kj}p_{kj} - x_{ij}| \leq u, -a \leq A_{kj} \leq a, 0 \leq C_{kj} \leq 2\}$, then according to (13), the cumulative distribution function (CDF) of $x'_{kj}$ can be derived as $G_{kj}(u) = P\{x'_{kj} \leq u\} = \iint_D \frac{1}{4a} dA_i dC_i$.

When $x_{ij} > 2p_{kj}$ and $u > a(x_{ij} - 2p_{kj})$, region $D$ is shown in Fig. 1(a) (marked by grey shading, similarly hereinafter). In this case, we have

$$G_{\alpha i}(u) = \frac{1}{4a}\left[2a + \frac{2u}{x_{ij}} + \int_{\frac{u}{x_{ij}}p_{kj}}^{a} \frac{u}{p_{kj}} \ln v \, dv - \left(\frac{x_{ij}}{p_{kj}} - 2\right)\left(a - \frac{u}{x_{ij}}\right)\right]$$
$$= \frac{1}{2} + \frac{1}{4a}\left[\frac{u}{p_{kj}} + \frac{u}{p_{kj}}\ln\frac{ax_{ij}}{u} - a\left(\frac{x_{ij}}{p_{kj}} - 2\right)\right]$$

When $x_{ij} > 2p_{kj}$ and $0 \leq u \leq a(x_{ij} - 2p_{kj})$, region $D$ is shown in Fig. 1(b). In this case, we have

$$G_{kj}(u) = \frac{1}{2} + \frac{u}{4ap_{kj}}\ln\frac{x_{ij}}{x_{ij} - 2p_{kj}}$$

When $p_{kj} \leq x_{ij} \leq 2p_{kj}$ and $u > a(2p_{kj} - x_{ij})$, region $D$ is shown in Fig. 1(c). In this case, we have

$$G_{kj}(u) = \frac{1}{2} + \frac{1}{4a}\left[\frac{u}{p_{kj}} + \frac{u}{p_{kj}}\ln\frac{ax_{ij}}{u} + a\left(2 - \frac{x_{ij}}{p_{kj}}\right)\right]$$

When $p_{kj} \leq x_i \leq 2p_{kj}$ and $0 \leq u \leq a(2p_{kj} - x_{ij})$, region $D$ is shown in Fig. 1(d). In this case, we have

$$G_{kj}(u) = \frac{1}{2} + \frac{1}{2a}\left[\frac{u}{p_{kj}} + \frac{u}{p_{kj}}\ln\frac{a\sqrt{x_{ij}(2p_{kj} - x_{ij})}}{u}\right]$$

Thus

$$g_{kj}(u) = \frac{dG_{kj}(u)}{du} = \begin{cases} \frac{1}{4ap_{kj}}\ln\frac{ax_{ij}}{u}, x_{ij} > 2p_{kj}, u > a(x_{ij} - 2p_{kj}) \\ \frac{1}{4ap_{kj}}\ln\frac{x_{kj}}{x_{kj} - 2p_{kj}}, x_{ij} > 2p_{kj}, 0 \leq u \leq a(x_{ij} - 2p_{kj}) \\ \frac{1}{4ap_{kj}}\ln\frac{ax_{ij}}{u}, p_{kj} \leq x_{ij} \leq 2p_{kj}, u > a(2p_{kj} - x_{ij}) \\ \frac{1}{2ap_{kj}}\ln\frac{a\sqrt{x_{ij}(2p_{kj} - x_{ij})}}{u}, p_{kj} \leq x_i \leq 2p_{kj}, 0 \leq u \leq a(2p_{kj} - x_{ij}) \end{cases}$$

Expanding the equations above to other cases where $x_{ij}, p_{kj}, u \in R$, and after a simple arrangement with absolute value sign, we finally get (15)-(16). □

Based on Theorem 1, we have

**Inference 1**: (1) The curve of $g_{kj}(u)$ is symmetrical about $u = p_{kj}$, and $g_{kj}(u)$ is non-decreasing within $(-\infty, p_{kj})$.

(2) $Ex'_{kj} = p_{kj}, Dx'_{kj} = \frac{a^2}{27}\left[\left(x_{ij}(t) - p_{kj}(t)\right)^2 + \frac{1}{3}p_{kj}^2(t)\right]$.

**Proof:** (1) This conclusion is obvious owing to the expression of $g_{kj}(u)$ (see (14)-(16)).

(2) The results can be obtained by calculating $\int_{-\infty}^{+\infty} u g_{kj}(u) du$ and $\int_{-\infty}^{+\infty} (u - p_{kj})^2 g_{kj}(u) du$. □

Based on (15) - (16), the curves of $g_{kj}(u)$ are shown in Fig. 2.

Furthermore, we can derive the PDF of $x_{ij}(t + 1)$ as:

**Theorem 2:** $h_{ij}(u) = (g_{1j} * g_{2j} * g_{3j})(3u)$

where the operator '$*$' denotes the convolution of two functions: $(f * g)(x) = \int_{-\infty}^{+\infty} f(u)g(x - u)du$.

**Proof**: the PDF of the sum of two independent random variables equals to the convolution of the PDFs of the two random variables [20]. □

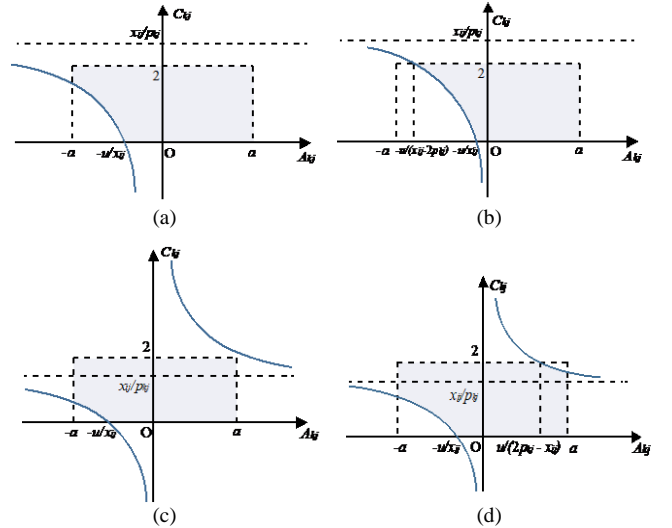

Fig. 1. Region $D$ (marked by grey shading), (a): $x_{ij} > 2p_{kj}$ and $u > a(x_{ij} - 2p_{kj})$; (b): $x_{ij} > 2p_{kj}$ and $0 \leq u \leq a(x_{ij} - 2p_{kj})$; (c): $p_{kj} \leq x_{ij} \leq 2p_{kj}$ and $u > a(2p_{kj} - x_{ij})$; (d): $p_{kj} \leq x_{ij} \leq 2p_{kj}$ and $0 \leq u \leq a(2p_{kj} - x_{ij})$.

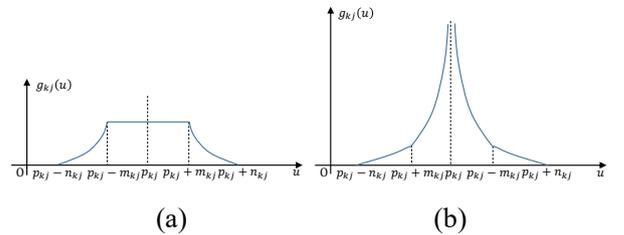

Fig. 2. The curves of $g_{kj}(u)$, (a): $m_{kj} > 0$, (b): $m_{kj} \leq 0$.

Based on Theorem 2 and Inference 1, we have:

**Inference 2**: (1) The curve of $h_{ij}(u)$ is symmetrical about $u = \frac{1}{3}\sum_{k=1}^{3} p_{kj}$, and $h_{ij}(u)$ is non-decreasing within $\left(-\infty, \frac{1}{3}\sum_{k=1}^{3} p_{kj}\right)$.

(2) $Ex_{ij}(t + 1) = \frac{1}{3}\sum_{k=1}^{3} p_{kj}, Dx_{ij}(t + 1) = \frac{a^2}{27}\sum_{k=1}^{3}\left[(x_{ij} - \right.$

$p_{kj})^2 + \frac{1}{3}p_{kj}^2\Big]$.

**Proof**: (1) The former part of the conclusion is obvious due to the characteristics of convolution. As for the latter part, we first proof $(g_{1j} * g_{2j})(u)$ is non-decreasing within $(-\infty, p_{1j} + p_{2j}]$.

For the sake of simplification, let $g_{10}(u) = g_{1j}(u + p_{1j})$, $g_{20}(u) = g_{2j}(u + p_{2j})$, then according to Inference 1(1), both $g_{10}(u), g_{20}(u)$ are non-decreasing on $R^-$, and we only need to prove $h_0(u) = (g_{10} * g_{20})(u)$ is non-decreasing on $R^-$. Let $u < u + \Delta u < 0$, we have

$h_0(u + \Delta u) - h_0(u)$
$= \int_{-\infty}^{+\infty} g_{10}(v)\big(g_{20}(u + \Delta u - v) - g_{20}(u - v)\big) dv$

Let $c(v) = g_{20}(u + \Delta u - v) - g_{20}(u - v)$, then obviously $c(2u + \Delta u - v) = -c(v)$, and when $v < u + \frac{1}{2}\Delta u$, $c(v) \leq 0$. Thus

$h_0(u + \Delta u) - h_0(u)$
$= \int_{-\infty}^{u+\frac{1}{2}\Delta u} g_{10}(v)c(v) dv + \int_{u+\frac{1}{2}\Delta u}^{+\infty} g_{10}(v)c(v) dv$
$= \int_{-\infty}^{u+\frac{1}{2}\Delta u} \big(g_{10}(v) - g_{10}(2u + \Delta u - v)\big)c(v) dv \geq 0$

i.e. $h_0(u)$ is non-decreasing on $R^-$, and consequently $(g_{1j} * g_{2j})(u)$ is non-decreasing within $(-\infty, p_{1j} + p_{2j})$. Similarly, $h_{ij}(u)$ is non-decreasing within $\left(-\infty, \frac{1}{3}\sum_{k=1}^{3} p_{kj}\right)$.

(2) According to Inference 1(1),

$\mathrm{E}x_{ij}(t+1) = \frac{1}{3}\sum_{k=1}^{3} \mathrm{E}x'_{kj} = \frac{1}{3}\sum_{k=1}^{3} p_{kj}$.

Since $x'_{kj}$ ($k = 1,2,3$) are independent with each other, according to Inference 1(2),

$\mathrm{D}x_{ij}(t+1) = \frac{1}{9}\sum_{k=1}^{3} \mathrm{D}x'_{kj} = \frac{a^2}{27}\sum_{k=1}^{3}\Big[(x_{ij} - p_{kj})^2 + \frac{1}{3}p_{kj}^2\Big]$.

□

### B. The Proposed BBGWO Based on Normal Distribution

Inference 2 shows that the PDF of $x_{ij}(t+1)$ is a unimodal function symmetrical about $u = \frac{1}{3}\sum_{k=1}^{3} p_{kj}$, whose curve is similar to that of the PDF of normal distribution. Thus, we can apply the bare bones mechanism to the original GWO and approximate the original solution updating equation (13) as:

$$x_{ij}(t+1) \sim N(\mu, \sigma^2) \quad (17)$$

$$\mu = \frac{1}{3}\sum_{k=1}^{3} p_{kj}, \sigma = \frac{a}{3\sqrt{3}}\sqrt{\sum_{k=1}^{3}\Big[(x_{ij} - p_{kj})^2 + \frac{1}{3}p_{kj}^2\Big]} \quad (18)$$

where the mean μ and standard deviation σ come from Inference 2 (2).

In this paper, the Monte Carlo method is applied to verify the validity of the approximation we proposed. Let $a = 2$, and $x_{ij}, p_{kj}$ are chosen randomly, the frequency histogram (group number = 80) of $x_{ij}(t+1)$ generated by Monte Carlo method (sample size = 100000) from (13) and the PDF of $x_{ij}(t+1)$ from (17) - (18) are shown in Fig. 3 (different subgraphs correspond to different values of $x_{ij}$ and $p_{kj}$).

It can be seen from Fig. 3 that the shape of the frequency histogram generated from (13) is similar to that of the PDF of $x_{ij}(t+1)$ from (17) - (18), and in most cases, they are highly similar to each other. This conclusion can be explained by central limit theorem, which means the properly normalized sum of independent random variables tends toward a normal distribution [20]. In this case, according to (13), $x_{ij}(t+1)$ is the sum of three independent random variables $x'_{kj}$, thus tends to a normal distribution.

Based on the discussions above, we can propose the BBGWO algorithm based on normal distribution by substituting the original solution updating equation of GWO (see (7) - (11)) by random variables generated by (17) - (18), i.e., applying the bare bones mechanism. The corresponding pseudocode of the proposed BBGWO is shown in Fig. 4, which is identical to that of the original GWO except for line 8 (shown in red), where in the original GWO, positions of current search agents are updated by (7) - (11).

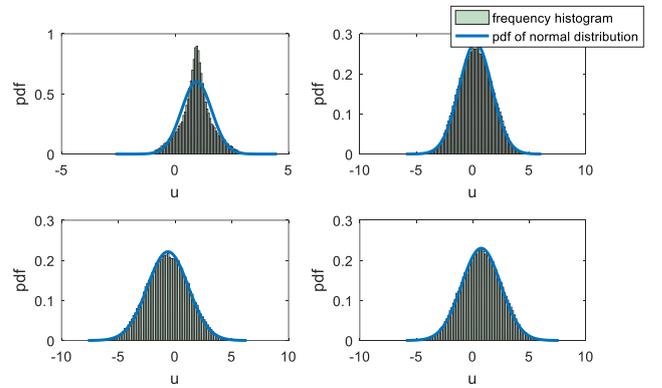

Fig. 3. Frequency histogram generated from (13) and PDF from (17) - (18).

---

1: *Initialize the grey wolf population* $\mathbf{x}_i (i = 1, ..., N)$
2: *Calculate the fitness of each search agent*
3: *Name the best search agent as* $\mathbf{p}_1$
4: *Name the second search agent as* $\mathbf{p}_2$
5: *Name the best search agent as* $\mathbf{p}_3$
6: **while** $t \leq T$
7:    **for** *each search agent*
8:      *Update the position of the current search agent by (12), (17)-(18)*
9:    **end for**
10:   *Calculate the fitness of all search agents*
11:   *Update* $\mathbf{p}_1, \mathbf{p}_2, \mathbf{p}_3$
12:   $t = t + 1$
13: **end while**
14: *return* $\mathbf{p}_1$ *and* $f(\mathbf{p}_1)$

Fig. 4. Pseudocode of the proposed BBGWO

## IV. RESULT AND DISCUSSION

In this paper, we set the parameters of GWO and BBGWO as: population size $N = 20$, maximum iterations $T = 500$, and step parameter $a$ updated by (12). All benchmark functions are given in Table I, which are all classical functions utilized by many researchers for testing the performance of optimization algorithms [21]-[23]. Each test function is optimized for 30 trials independently, and the average value and variance of the 30 simulation results from the GWO and BBGWO are calculated, aiming at measuring the convergence speed and stability of an algorithm. Moreover, for each benchmark function, the times that an algorithm finds the global optimum (i.e., the difference between the result from the algorithm and the minimum value of the benchmark function is less than 0.001) is calculated as well. All simulation results are shown in Table II.

It can be seen from Table II that for most cases the performance of the BBGWO is similar to that of the GWO. Specifically speaking, for those cases in which the GWO is able to find the global optimum, the BBGWO is able to find the global optimum as well, and both the convergence speed and variance of these two algorithms are similar to each other. As for cases in which the GWO could hardly find the global optimum (for instance, the functions 4, 6, 10 and 12), the BBGWO could hardly find the global optimum either. Thus, the theoretical value and practical significance of the proposed BBGWO are verified.

## V. CONCLUSION

In this paper, we apply the bare bones mechanism to the original GWO to obtain a novel algorithm called BBGWO, in which the complex solution updating mechanism is replaced by a random vector that obeys a normal distribution. The theoretical analysis and simulation results show that the proposed BBGWO algorithm has the same optimization performance as the the original GWO, which means the proposed BBGWO can be applied to solve engineering optimization problems independently, and on the other hand, the BBGWO may serve as a fundamental tool in studying the mechanism of GWO.

TABLE I
BENCHMARK FUNCTIONS

| No. | Minimum value | Average value | | Variance | | Times to find the global optimum | |
|---|---|---|---|---|---|---|---|
| | | GWO | BBGWO | GWO | BBGWO | GWO | BBGWO |
| 1 | 0 | $5.0976 \times 10^{-23}$ | $3.2750 \times 10^{-23}$ | $6.3039 \times 10^{-45}$ | $6.1933 \times 10^{-45}$ | 30 | 30 |
| 2 | 0 | $3.5800 \times 10^{-14}$ | $2.3377 \times 10^{-14}$ | $6.5943 \times 10^{-28}$ | $2.8613 \times 10^{-28}$ | 30 | 30 |
| 3 | 0 | 0.0011 | $7.3808 \times 10^{-4}$ | $7.7187 \times 10^{-6}$ | $4.9808 \times 10^{-6}$ | 23 | 27 |
| 4 | 0 | 28.1752 | 28.1149 | 0.3067 | 0.3458 | 0 | 0 |
| 5 | 0 | 0.8667 | 1.0333 | 1.4299 | 1.2747 | 16 | 12 |
| 6 | -7286.2 | -5841.2 | -5897.6 | 68435 | 92390 | 0 | 0 |
| 7 | 0 | 5.9172 | 4.4231 | 24.2071 | 29.7328 | 7 | 11 |
| 8 | 0 | $8.5111 \times 10^{-13}$ | $6.0325 \times 10^{-13}$ | $2.1907 \times 10^{-25}$ | $1.5212 \times 10^{-25}$ | 30 | 30 |
| 9 | 0 | 0.0059 | 0.0062 | $1.0906 \times 10^{-4}$ | $1.1700 \times 10^{-4}$ | 22 | 22 |
| 10 | -29.6248 | -14.6104 | -12.4215 | 13.1229 | 10.8998 | 0 | 0 |
| 11 | -1.0316 | -1.0316 | -1.0316 | $2.4452 \times 10^{-15}$ | $2.2518 \times 10^{-15}$ | 30 | 30 |
| 12 | 0 | 0.0854 | 0.0879 | 0.0109 | 0.0096 | 2 | 1 |

TABLE II
SIMULATION RESULTS

| No. | Expression | $D$ | Domain | Minimum value | Any local optimum |
|---|---|---|---|---|---|
| 1 | $f(x) = \sum_{j=1}^{D} x_j^2$ | 30 | $[-100,100]^D$ | 0 | No |
| 2 | $f(x) = \sum_{j=1}^{D}|x_j| + |\prod_{j=1}^{D} x_j|$ | 30 | $[-10,10]^D$ | 0 | No |
| 3 | $f(x) = \sum_{j=1}^{D}\left(\sum_{k=1}^{j} x_k\right)^2$ | 30 | $[-100,100]^D$ | 0 | No |
| 4 | $f(x) = \sum_{j=1}^{D-1}(x_j - x_{j+1}^2)^2 + \sum_{j=1}^{D}(x_j - 1)^2$ | 30 | $[-30,30]^D$ | 0 | No |
| 5 | $f(x) = \sum_{j=1}^{D}\lfloor x_j \rfloor^2$ | 30 | $[-100,100]^D$ | 0 | No |
| 6 | $f(x) = -\sum_{j=1}^{D} x_j \sin\sqrt{|x_j|}$ | 30 | $[-500,500]^D$ | -7286.2 | Yes |
| 7 | $f(x) = \sum_{j=1}^{D}(x_j^2 - 10\cos 2\pi x_j + 10)$ | 30 | $[-10,10]^D$ | 0 | Yes |
| 8 | $f(x) = 20 + e - 20e^{-0.2\sqrt{\frac{1}{D}\sum_{j=1}^{D} x_j^2}} - e^{\frac{1}{D}\sum_{j=1}^{D}\cos 2\pi x_j}$ | 30 | $[-20,20]^D$ | 0 | Yes |
| 9 | $f(x) = 1 + \frac{1}{4000}\sum_{j=1}^{D} x_j^2 - \prod_{j=1}^{D}\cos\frac{x_j}{\sqrt{j}}$ | 30 | $[-600,600]^D$ | 0 | Yes |
| 10 | $f(x) = -\sum_{j=1}^{D}\sin x_j \left(\sin\frac{jx_j^2}{\pi}\right)^{20}$ | 30 | $[0,\pi]^D$ | -29.6248 | Yes |
| 11 | $f(x) = 4x_1^2 - 2.1x_1^4 + \frac{1}{3}x_1^6 + x_1 x_2 - 4x_2^2 + 4x_2^4$ | 2 | $[-5,5]^D$ | -1.0316 | Yes |
| 12 | $f(x) = \sum_{i=1}^{n}\left[q_i - x_1 e^{-\frac{(p_i - x_3)^2}{2x_2^2}} - x_4 e^{-\frac{(p_i - x_6)^2}{2x_5^2}}\right]^2$ | 6 | $[0,10]^D$ | 0 | Yes |